\documentclass[12pt]{amsart}
\usepackage{graphicx,dbnsymb,amssymb,picins,floatflt,color,epic,eepic}
\usepackage[all]{xy}

\theoremstyle{plain}

\newtheorem{proposition}{Proposition}[section]

\newtheorem{lemma}[proposition]{Lemma}

\theoremstyle{definition}

\theoremstyle{remark}

\newcommand{\mathmode}[1]{$#1$}
\newlength{\standardunitlength}
\catcode`\@=11
\long\def\@makecaption#1#2{%
    \vskip 10pt
    \setbox\@tempboxa\hbox{%\ifvoid\tinybox\else\box\tinybox\fi
      \small\sf{\bfcaptionfont #1. }\ignorespaces #2}%
    \ifdim \wd\@tempboxa >\captionwidth {%
        \rightskip=\@captionmargin\leftskip=\@captionmargin
        \unhbox\@tempboxa\par}%
      \else
        \hbox to\hsize{\hfil\box\@tempboxa\hfil}%
    \fi}
\font\bfcaptionfont=cmssbx10 scaled \magstephalf
\newdimen\@captionmargin\@captionmargin=2\parindent
\newdimen\captionwidth\captionwidth=\hsize
\catcode`\@=12

\def\qed{{\hfill\text{$\Box$}}}

\newlength{\globalparindent}
\def\llbracket{\left[\!\!\left[}
\def\rrbracket{\right]\!\!\right]}

\def\bbQ{{\mathbb Q}}

\def\bbZ{{\mathbb Z}}

\def\calC{{\mathcal C}}

\def\calF{{\mathcal F}}
\def\calG{{\mathcal G}}

\newcommand{\Cobdl}{{\mathcal Cob}_{\bullet/l}}

\newcommand{\Cobl}{{\mathcal Cob}_{/l}}

\newcommand{\Kh}{{\text{\it Kh}}}

\newcommand{\Kom}{\operatorname{Kom}}

\newcommand{\Mat}{\operatorname{Mat}}

\renewcommand{\qed}{~\hfill$\square$}

\begin{document}
\newdimen\captionwidth\captionwidth=\hsize

\title{Fast Khovanov Homology Computations}

\author{Dror~Bar-Natan}
\address{
  Department of Mathematics\\
  University of Toronto\\
  Toronto Ontario M5S 3G3\\
  Canada
}
\email{drorbn@math.toronto.edu}
\urladdr{http://www.math.toronto.edu/\~{}drorbn}

\date{
  First edition: Jun.{} 13, 2006.
  This edition: Jun.~13,~2006
}


\thanks{This work was partially supported by NSERC grant RGPIN 262178.
  Electronic version: {\tt
  http://\linebreak[0]www.math.toronto.edu/\linebreak[0]$\sim$drorbn/\linebreak[0]papers/\linebreak[0]FastKh/} and arXiv:math.GT/0606318.
}

\begin{abstract}
  We introduce a {\em local} algorithm for Khovanov Homology
  computations --- that is, we explain how it is possible to ``cancel''
  terms in the Khovanov complex associated with a (``local'') tangle,
  hence canceling the many associated ``global'' terms in one swoosh
  early on. This leads to a dramatic improvement in computational
  efficiency.  Thus our program can rapidly compute certain Khovanov
  homology groups that otherwise would have taken centuries to evaluate.
\end{abstract}

\dedicatory{To Lou Kauffman, who gave us $\slashoverback\mapsto
A\,\smoothing+A^{-1}\hsmoothing$.}

\maketitle

\tableofcontents

\section{Introduction} \label{sec:intro}

%\begin{floatingfigure}[r]{60mm}
%\hspace{-5mm}\includegraphics[width=60mm]{figs/SusanWilliamsShirtPin_640.ps}
%\caption{A knot cut in two~\cite{Bar-Natan:ShirtPin}.} \label{fig:ShirtPin}
%\end{floatingfigure}

\parpic[r]{$\begin{array}{r}
  \includegraphics[width=45mm]{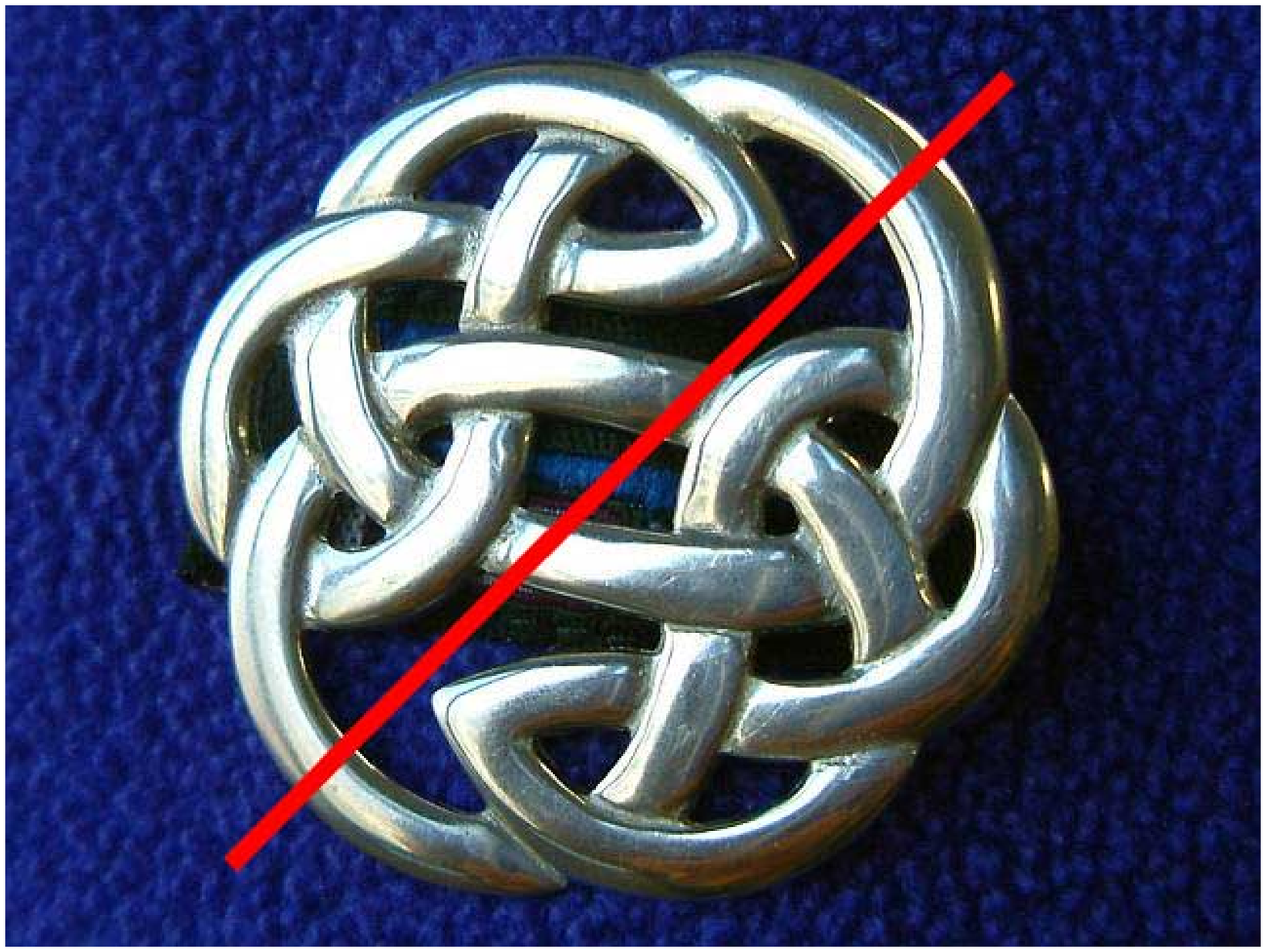} \\
  \text{\tiny Image source:~\cite{Bar-Natan:ShirtPin}.}
\end{array}$}
The ``divide and conquer'' approach to computation, as applied to
knot theory, goes roughly as follows. Suppose a certain knot invariant
takes an exponential amount of time to compute, so on a knot with $n$
crossings ($n=14$ on the right), it takes roughly $C_1^n$ operations,
where $C_1$ is some constant. Suppose also that that same knot invariant
can also be computed ``in halves''; i.e., it makes sense to compute a
``partial'' invariant of the left half of the knot, consisting of just
$n/2$ crossings, and likewise for the right half of the knot. Then, with
some luck, it takes just $C_1^{n/2}$ operations to do each half and if
the assembly of the half-computations into the full one is cheap then
the whole computation takes $2C_1^{n/2}$ operations, a lot less than the
original $C_1^n$ operations.  Of course, by iterating this procedure
one can save even more and carry out the computation in $4C_1^{n/4}$
operations, or even just $8C_1^{n/8}$ operations. At the limit the
computation time becomes linear linear in $n$, at least if one ignores
the costs of cutting and of assembly.

In reality, cutting is indeed cheap but assembly isn't. Often the
invariant of each half knot has to be quite complicated in order to
allow for its pairing with every conceivable ``other half''. This often
means that each half-knot invariant must take value in a space whose
dimension grows exponentially in the number of strands $b$ that connect the
two halves ($b=4$ in Figure~\ref{fig:ShirtPin}). Thus each half
computation takes at least $C_2^b$ operations and if luck strikes, it
doesn't take much longer. Typically we can expect $b$ to be roughly the
``width'' of the knot and since we are in the plane, we can expect $b$
to be around $\sqrt{n}$. Thus realistically divide and conquer may
reduce $C_1^n$ to $C_2^{\sqrt{n}}$. The latter is still very big, but
it is a lot smaller than the former.

In fact, the advantage of ``divide and conquer'' is so big that it is
worthwhile to try this approach even if the assembly cost is more than
$C_2^b$ or even if no good estimates for the assembly cost at all exist, as
is the case at hand in this paper.

This paper applies the ``divide and conquer'' approach to the
computation of Khovanov homology~\cite{Khovanov:Categorification,
Bar-Natan:Categorification}. We start in Section~\ref{sec:QuickReview}
with a quick review of the local Khovanov theory
of~\cite{Bar-Natan:Cobordisms}, which amounts to a definition of
``Khovanov homology'' for half-knots (i.e., for tangles), along with
the ``horizontal composition'' technique necessary for the assembly of
the invariants of two tangles into the invariant of their side-by-side
composition.

\begin{figure}
\[ \includegraphics[height=1.0in]{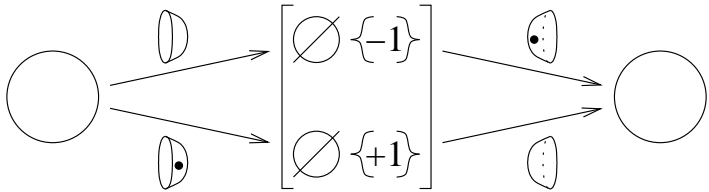} \]
\caption{Delooping.} \label{fig:Delooping}
\end{figure}

In Section~\ref{sec:Tools} we introduce two simple tools, delooping and
Gaussian elimination, that allow us to ``simplify'' the invariants of
tangles. These tools are the keys to the whole paper, as they reduce
the complexity of the Khovanov complex associated with a tangle and
thus allow for much easier horizontal composition assembly.  For the
impatient reader, delooping and Gaussian elimination are depicted in
Figures~\ref{fig:Delooping} and~\ref{fig:GaussianElimination} here. If
you understand these sketches you've understood the whole paper.

\begin{figure}
\[ \begin{array}{c}
  \xymatrix@C=2cm{
    \left[C\right]
    \ar[r]^{\begin{pmatrix}\alpha \\ \beta\end{pmatrix}} &
    {\begin{bmatrix}b_1 \\ D\end{bmatrix}}
    \ar[r]^{\begin{pmatrix} \phi & \delta \\ \gamma & \epsilon \end{pmatrix}} &
    {\begin{bmatrix}b_2 \\ E\end{bmatrix}}
    \ar[r]^{\begin{pmatrix} \mu & \nu \end{pmatrix}} &
    \left[F\right]
  } \\
  \text{is isomorphic to the (direct sum) complex} \\
  \xymatrix@C=3cm{
    \left[C\right]
    \ar[r]^{\begin{pmatrix}0 \\ \beta\end{pmatrix}} &
    {\begin{bmatrix}b_1 \\ D\end{bmatrix}}
    \ar[r]^{\begin{pmatrix}
      \phi & 0 \\ 0 & \epsilon-\gamma\phi^{-1}\delta
    \end{pmatrix}} &
    {\begin{bmatrix}b_2 \\ E\end{bmatrix}}
    \ar[r]^{\begin{pmatrix} 0 & \nu \end{pmatrix}} &
    \left[F\right]
  }
\end{array} \]
\caption{Gaussian elimination.} \label{fig:GaussianElimination}
\end{figure}

In Section~\ref{sec:Algorithm} we describe our algorithm. As an
illustration, in Section~\ref{sec:FigureEight} we ``run'' the algorithm
on the figure eight knot. Then in Section~\ref{sec:Faster} we describe
in even faster variant of the algorithm and in
Section~\ref{sec:ComputerPrograms} we mention the two available
implementations and exhibit a sample computation.  The final
Section~\ref{sec:Reidemeister} quickly explains how the tools used in
this paper lead to a completely automated proof of the invariance of
Khovanov homology under Reidemeister moves.

\section{History and Acknowledgement} When I did the research for my
paper~\cite{Bar-Natan:Cobordisms} part of my motivation was to set up
the framework allowing for ``divide and conquer'' computations as
described in Section~\ref{sec:intro}. But by the time I finished
writing~\cite{Bar-Natan:Cobordisms} I completely forgot about that part
of my motivation. I wish to thank Marco Mackaay for reminding me!
Further thanks to Jeremy Green for implementing the algorithm described
here and to Louis Leung and Gad Naot for some comments and suggestions.

\section{A quick review of the local Khovanov theory}
\label{sec:QuickReview}

\parpic[r]{$\begin{array}{c}
  \begin{picture}(0,0)%
\includegraphics{figs/Main2.pstex}%
\end{picture}%
%
%  pstex_opts: -m 0.65 
%
\setlength{\unitlength}{2565sp}%
\begingroup\makeatletter\ifx\SetFigFont\undefined%
\gdef\SetFigFont#1#2#3#4#5{%
  \reset@font\fontsize{#1}{#2pt}%
  \fontfamily{#3}\fontseries{#4}\fontshape{#5}%
  \selectfont}%
\fi\endgroup%
\begin{picture}(5894,4144)(1189,-3968)
\put(1801,-1186){\makebox(0,0)[b]{\smash{\SetFigFont{8}{9.6}{\rmdefault}{\mddefault}{\updefault}{\color[rgb]{0,0,0}$(n_+,n_-)=(2,0)$}%
}}}
\end{picture}

\end{array}$}
Let us briefly recall the definition of the Khovanov complex for tangles,
following \cite{Bar-Natan:Cobordisms}. Given an $n$-crossing tangle $T$
with boundary $\partial T$, such as the 2-crossing tangle displayed on the
right, one constructs an $n$-dimensional ``cube'' of smoothings and
cobordisms between them (as illustrated on the right). This cube is
then ``flattened'' to a ``formal complex'' $\llbracket T\rrbracket$ in
the additive category $\Cobdl^3(\partial T)$ whose objects are
formally graded smoothings with boundary $\partial T$ and whose
morphisms are formal linear combinations of cobordisms whose tops and
bottoms are smoothings and whose side boundaries are $I\times\partial
T$, modulo some local relations.

The Khovanov complex $\Kh(T)$ of $T$ is obtained from $\llbracket
T\rrbracket$ by some minor further degree and height shifts depending only
of the numbers $n_\pm$ of over- and under-crossings in $T$
(see~\cite[Definition~6.4]{Bar-Natan:Cobordisms}). It is a member
in the category $\Kom(\Mat(\Cobdl^3(\partial T)))$ of complexes of formal
direct sums of objects in $\Cobdl^3(\partial T)$ and it is invariant up to
(formal) homotopies.

For simplicity we are using as the basis to our story one of the
simpler cobordism categories $\Cobdl^3$ that appear
in~\cite{Bar-Natan:Cobordisms}, rather than the most general one,
$\Cobl^3$ (a fuller treatment appears in~\cite{Naot}). It is worthwhile
to repeat here the local relations that appear in the definition of
$\Cobdl^3$ (see~\cite[Section 11.2]{Bar-Natan:Cobordisms}):

\begin{equation} \label{eq:LocalRelations}
\begin{array}{c}
  \begin{array}{c}
    \includegraphics[height=1cm]{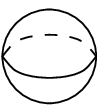}
  \end{array}\hspace{-2mm}=0,
  \qquad\qquad
  \begin{array}{c}
    \includegraphics[height=1cm]{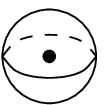}
  \end{array}\hspace{-2mm}=1,
  \qquad\qquad
  \begin{array}{c}\includegraphics[height=10mm]{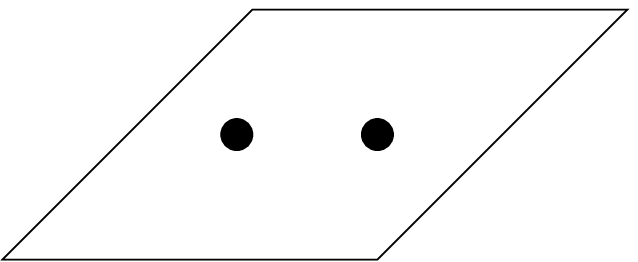}\end{array}
  \hspace{-4mm}=0,
\\
  \text{and}\qquad
  \begin{array}{c}\includegraphics[height=10mm]{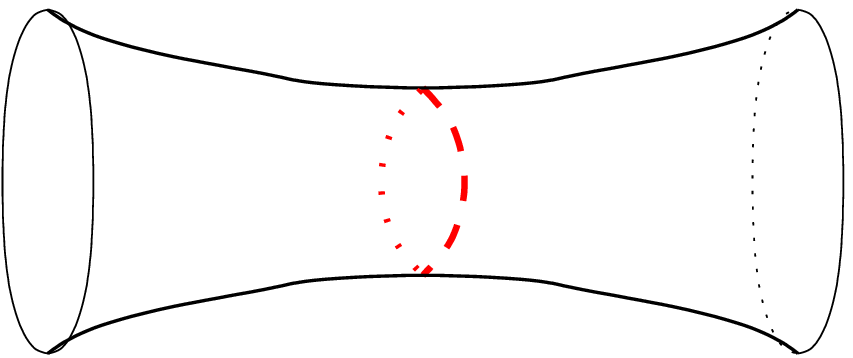}\end{array}
  =\begin{array}{c}\includegraphics[height=10mm]{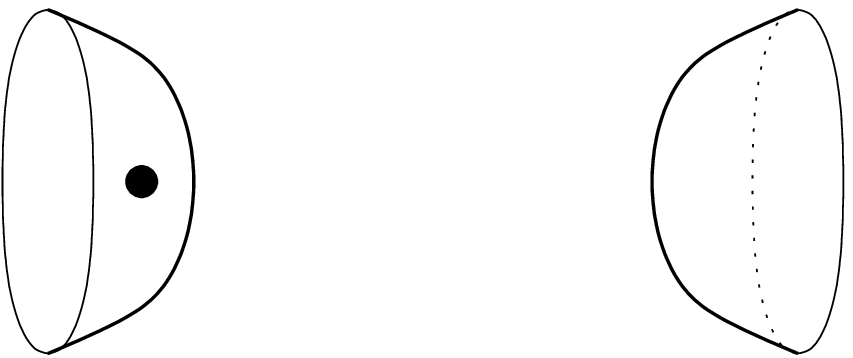}\end{array}
  +\begin{array}{c}\includegraphics[height=10mm]{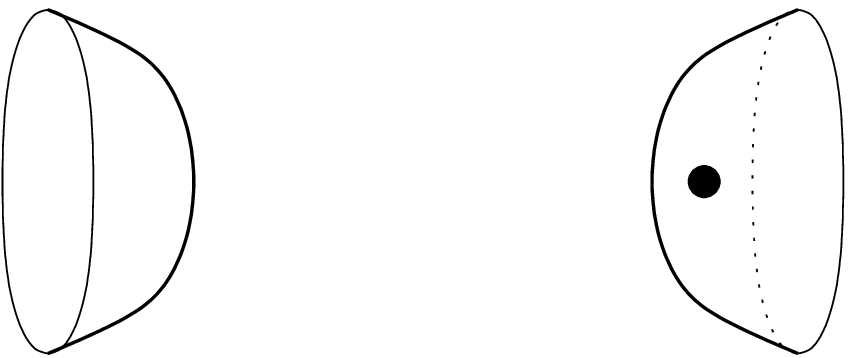}\end{array}.
\end{array}
\end{equation}

Also recall from~\cite[Section 5]{Bar-Natan:Cobordisms} that
$\llbracket\cdot\rrbracket$, and hence $\Kh(\cdot)$, are planar algebra
morphisms. That is, if $T_1$ and $T_2$ are tangles and $D(T_1,T_2)$
denotes one of their side-by-side compositions (a side by side
placement of $T_1$ and $T_2$ while joining some of their ends in a
certain way prescribed by a planar arc diagram $D$), then $\llbracket
D(T_1,T_2)\rrbracket=D(\llbracket T_1\rrbracket,\llbracket
T_2\rrbracket)$. Here as in~\cite[Section 5]{Bar-Natan:Cobordisms}
$D(\llbracket T_1\rrbracket,\llbracket T_2\rrbracket)$ is the ``tensor
product'' operation induced on formal complexes by the horizontal
composition operation $D$ on the canopoly $\Cobdl^3$. In exactly the same
sense we also have that $\Kh(D(T_1,T_2))=D(\Kh(T_1),\Kh(T_2))$.

Thus Khovanov homology is ready for a divide and conquer computation.
It makes sense for half-knots (tangles) and there is a composition rule
that takes the invariants of the halves and produces the invariant of
the whole.  But as it stands there is no advantage (yet) for this
computation method.  If $T_i$ (for $i=1,2$) has $n_i$ crossings, the
Khovanov cube for $T_i$ consists of $2^{n_i}$ vertices, thus $\llbracket
T_i\rrbracket$ involves $2^{n_i}$ objects and thus $D(\llbracket
T_1\rrbracket,\llbracket T_2\rrbracket)$ involves $2^{n_1}\cdot
2^{n_2}=2^{n_1+n_2}$ objects, exactly as many as in $\llbracket
D(T_1,T_2)\rrbracket$, and nothing has been saved.

\section{The tools: delooping and Gaussian elimination} \label{sec:Tools}

To overcome the difficulty from the previous paragraph we need to learn
how to simplify $\llbracket T_1\rrbracket$ and $\llbracket T_2\rrbracket$
(modulo homotopy) {\em before} taking their tensor product and thus before
the biggest number ($2^{n_1+n_2}$) is encountered. For this we need tools
for simplifying complexes over the category $\Cobdl^3$. These tools are
Lemma~\ref{lem:Delooping} and Lemma~\ref{lem:GaussianElimination} below.

\begin{lemma} \label{lem:Delooping} (Delooping) If an object $S$ in $\Cobdl^3$
contains a closed loop $\ell$, then it is isomorphic (in
$\Mat(\Cobdl^3)$) to the direct sum of two copies $S'\{+1\}$ and $S'\{-1\}$
of $S$ in which $\ell$ is removed, one taken with a degree shift of
$+1$ and one with a degree shift of $-1$. Symbolically, this reads
$\bigcirc\equiv\emptyset\{+1\}\oplus\emptyset\{-1\}$.
\end{lemma}

\begin{proof} Here are the isomorphisms:
\begin{equation} \label{eq:Delooping}
  \begin{array}{c}
    \includegraphics[height=1.0in]{figs/Delooping.eps}
  \end{array}
\end{equation}
It is easy to verify using (all!) the relations
in~\eqref{eq:LocalRelations} that the two possible compositions of the
morphisms above are both equal to the identity morphisms of the
relevant objects.  (This is the only place in this paper where the
relations in~\eqref{eq:LocalRelations} are used). \qed
\end{proof}

\begin{lemma} \label{lem:GaussianElimination} (Gaussian elimination, made
abstract) If $\phi:b_1\to b_2$ is an isomorphism (in some additive
category $\calC$), then the four term complex segment in $\Mat(\calC)$
\begin{equation} \label{eq:BeforeGE}
  \xymatrix@C=2cm{
    \cdots\ 
    \left[C\right]
    \ar[r]^{\begin{pmatrix}\alpha \\ \beta\end{pmatrix}} &
    {\begin{bmatrix}b_1 \\ D\end{bmatrix}}
    \ar[r]^{\begin{pmatrix}
      \phi & \delta \\ \gamma & \epsilon
    \end{pmatrix}} &
    {\begin{bmatrix}b_2 \\ E\end{bmatrix}}
    \ar[r]^{\begin{pmatrix} \mu & \nu \end{pmatrix}} &
    \left[F\right] \  \cdots
  }
\end{equation}
is isomorphic to the (direct sum) complex segment
\begin{equation} \label{eq:DuringGE}
  \xymatrix@C=3cm{
    \cdots\ 
    \left[C\right]
    \ar[r]^{\begin{pmatrix}0 \\ \beta\end{pmatrix}} &
    {\begin{bmatrix}b_1 \\ D\end{bmatrix}}
    \ar[r]^{\begin{pmatrix}
      \phi & 0 \\ 0 & \epsilon-\gamma\phi^{-1}\delta
    \end{pmatrix}} &
    {\begin{bmatrix}b_2 \\ E\end{bmatrix}}
    \ar[r]^{\begin{pmatrix} 0 & \nu \end{pmatrix}} &
    \left[F\right] \  \cdots
  }.
\end{equation}
Both these complexes are homotopy equivalent to the (simpler) complex
segment
\begin{equation} \label{eq:AfterGE}
  \xymatrix@C=3cm{
    \cdots\ 
    \left[C\right]
    \ar[r]^{\left(\beta\right)} &
    {\left[D\right]}
    \ar[r]^{\left(\epsilon-\gamma\phi^{-1}\delta\right)} &
    {\left[E\right]}
    \ar[r]^{\left(\nu\right)} &
    \left[F\right] \  \cdots
  }.
\end{equation}
Here $C$, $D$, $E$ and $F$ are arbitrary columns of objects in $\calC$
and all Greek letters (other than $\phi$) represent arbitrary matrices
of morphisms in $\calC$ (having the appropriate dimensions, domains
and ranges); all matrices appearing in these complexes are block-matrices
with blocks as specified. $b_1$ and $b_2$ are billed here as individual
objects of $\calC$, but they can equally well be taken to be columns of
objects provided (the morphism matrix) $\phi$ remains invertible.
\end{lemma}

\begin{proof} The two $2\times 2$ (block) matrices in~\eqref{eq:BeforeGE}
and~\eqref{eq:DuringGE} differ by invertible row and column operations
(i.e., by a ``change of basis''). When the corresponding column and row
operations are performed on $\begin{pmatrix} \mu & \nu \end{pmatrix}$
and on $\begin{pmatrix}\alpha \\ \beta\end{pmatrix}$, the results are
$\begin{pmatrix}\mu-\nu\gamma\phi^{-1}&\nu\end{pmatrix}=\begin{pmatrix}
0 & \nu \end{pmatrix}$ and $\begin{pmatrix}\alpha-\phi^{-1}\delta\beta\\
\beta\end{pmatrix}=\begin{pmatrix} 0 \\ \beta\end{pmatrix}$ respectively
(note that $\mu\phi-\nu\gamma=0$ and $\phi\alpha-\delta\beta=0$ as
in the original complex the differential squares to $0$). Hence the complexes
in~\eqref{eq:BeforeGE} and~\eqref{eq:DuringGE} are isomorphic.

The complexes in~\eqref{eq:DuringGE} and~\eqref{eq:AfterGE} differ by the
removal of a contractible direct summand
$\xymatrix{0\ar[r]&b_1\ar[r]^\phi&b_2\ar[r]&0}$ (remember that $\phi$ is
invertible). Hence they are homotopy equivalent. \qed
\end{proof}

\section{The algorithm} \label{sec:Algorithm}

We thus have a procedure for simplifying complexes $\Omega$ with objects in
$\Cobdl$ (and thus, along with ``divide and conquer'', we have a
potentially fast way of computing Khovanov homologies):

\begin{itemize}
\item Whenever an object in $\Omega$ contains a closed loop, double it up
using Lemma~\ref{lem:Delooping}, removing the closed loop and inserting
$\pm 1$ degree shifts. This done, $\Omega$ becomes $\Omega'$.
\item $\Omega'$ is ``bigger'' than $\Omega$, but it is made up of many
fewer possible objects. Thus it is likely that many morphisms in $\Omega'$
are isomorphisms. Whenever you find one, cancel it using
Lemma~\ref{lem:GaussianElimination}. Call the result of doing this
iteratively $\Omega''$.
\end{itemize}

A priori, we cannot guarantee that $\Omega''$ will be simpler than
$\Omega$. But experimentation shows that it is, as seen in the next few
sections.

\section{The figure eight knot} \label{sec:FigureEight}

\parpic[r]{$
  \setlength{\unitlength}{0.75\standardunitlength}
  \begin{array}{c}
    {\begingroup\makeatletter\ifx\SetFigFont\undefined%
\gdef\SetFigFont#1#2#3#4#5{%
  \reset@font\fontsize{#1}{#2pt}%
  \fontfamily{#3}\fontseries{#4}\fontshape{#5}%
  \selectfont}%
\fi\endgroup%
{\renewcommand{\dashlinestretch}{30}
\begin{picture}(1554,1602)(0,-10)
\put(481,975){\makebox(0,0)[b]{\smash{{\mathmode{-}}}}}
\put(1456,600){\makebox(0,0)[b]{\smash{{\mathmode{+}}}}}
\put(1456,900){\makebox(0,0)[b]{\smash{{\mathmode{+}}}}}
\dashline{60.000}(931,1575)(931,75)
\path(31,975)(331,675)
\path(31,675)(106,750)
\path(331,975)(256,900)
\path(319.640,1006.066)(256.000,900.000)(362.066,963.640)
\path(331,975)(631,675)
\path(331,675)(406,750)
\path(342.360,643.934)(406.000,750.000)(299.934,686.360)
\path(631,975)(556,900)
\put(181,975){\makebox(0,0)[b]{\smash{{\mathmode{-}}}}}
\path(1081,825)(1381,525)
\put(1231,0){\makebox(0,0)[lb]{\smash{{\mathmode{T_R}}}}}
\path(1081,525)(1156,600)
\path(1381,825)(1306,750)
\path(1369.640,856.066)(1306.000,750.000)(1412.066,813.640)
\path(1081,1125)(1381,825)
\path(1081,825)(1156,900)
\path(1092.360,793.934)(1156.000,900.000)(1049.934,836.360)
\path(1381,1125)(1306,1050)
\path(1081,1125)(1080,1126)(1076,1129)
	(1066,1137)(1051,1150)(1032,1164)
	(1011,1177)(991,1189)(973,1198)
	(955,1203)(939,1205)(923,1204)
	(906,1200)(894,1196)(882,1190)
	(869,1182)(854,1173)(839,1161)
	(821,1147)(802,1131)(780,1113)
	(757,1093)(734,1071)(710,1049)
	(687,1028)(667,1010)(652,995)
	(640,984)(634,978)(631,975)
\path(1081,525)(1080,524)(1076,521)
	(1066,513)(1051,500)(1032,486)
	(1011,473)(991,461)(973,452)
	(955,447)(939,445)(923,446)
	(906,450)(894,454)(882,460)
	(869,468)(854,477)(839,489)
	(821,503)(802,519)(780,537)
	(757,557)(734,579)(710,601)
	(687,622)(667,640)(652,655)
	(640,666)(634,672)(631,675)
\path(1381,525)(1383,523)(1386,518)
	(1392,510)(1400,500)(1408,487)
	(1417,472)(1425,457)(1431,443)
	(1435,430)(1437,417)(1436,405)
	(1432,394)(1426,384)(1418,373)
	(1406,362)(1396,355)(1384,347)
	(1371,338)(1357,330)(1341,321)
	(1323,312)(1304,303)(1284,294)
	(1262,284)(1239,275)(1215,266)
	(1191,257)(1166,248)(1141,240)
	(1115,232)(1090,224)(1064,217)
	(1040,210)(1015,204)(991,198)
	(967,192)(944,187)(920,183)
	(897,179)(873,175)(849,172)
	(825,169)(800,167)(775,165)
	(749,164)(723,163)(698,163)
	(672,164)(646,165)(621,167)
	(596,170)(572,173)(548,177)
	(525,181)(503,186)(481,192)
	(460,198)(439,205)(419,212)
	(396,222)(374,231)(351,242)
	(329,254)(306,267)(284,280)
	(261,295)(238,310)(216,325)
	(195,341)(174,357)(154,373)
	(136,389)(119,405)(103,420)
	(88,435)(75,449)(63,462)
	(53,475)(44,487)(32,504)
	(24,520)(18,536)(14,552)
	(12,568)(12,585)(13,603)
	(16,621)(20,639)(24,654)
	(28,665)(30,672)(31,675)
\path(1381,1125)(1383,1127)(1386,1132)
	(1392,1140)(1400,1150)(1408,1163)
	(1417,1178)(1425,1193)(1431,1207)
	(1435,1220)(1437,1233)(1436,1245)
	(1432,1256)(1426,1266)(1418,1277)
	(1406,1287)(1396,1295)(1384,1303)
	(1371,1312)(1357,1320)(1341,1329)
	(1323,1338)(1304,1347)(1284,1356)
	(1262,1366)(1239,1375)(1215,1384)
	(1191,1393)(1166,1402)(1141,1410)
	(1115,1418)(1090,1426)(1064,1433)
	(1040,1440)(1015,1446)(991,1452)
	(967,1458)(944,1462)(920,1467)
	(897,1471)(873,1475)(849,1478)
	(825,1481)(800,1483)(775,1485)
	(749,1486)(723,1487)(698,1487)
	(672,1486)(646,1485)(621,1483)
	(596,1480)(572,1477)(548,1473)
	(525,1469)(503,1464)(481,1458)
	(460,1452)(439,1445)(419,1437)
	(396,1428)(374,1419)(351,1408)
	(329,1396)(306,1383)(284,1370)
	(261,1355)(238,1340)(216,1325)
	(195,1309)(174,1293)(154,1277)
	(136,1261)(119,1245)(103,1230)
	(88,1215)(75,1201)(63,1188)
	(53,1175)(44,1162)(32,1146)
	(24,1130)(18,1114)(14,1098)
	(12,1082)(12,1065)(13,1047)
	(16,1029)(20,1011)(24,996)
	(28,985)(30,978)(31,975)
\put(31,0){\makebox(0,0)[lb]{\smash{{\mathmode{T_L}}}}}
\end{picture}
} } 
  \end{array}
$}
Our first example is the figure eight knot, cut in half into two
tangles $T_L$ and $T_R$ as shown on the right. $T_L$ is the tangle
$
  \setlength{\unitlength}{0.5\standardunitlength}
  \begin{array}{c}  \hspace{-2.8mm}
    \raisebox{-2pt}{%
\begingroup\makeatletter\ifx\SetFigFont\undefined%
\gdef\SetFigFont#1#2#3#4#5{%
  \reset@font\fontsize{#1}{#2pt}%
  \fontfamily{#3}\fontseries{#4}\fontshape{#5}%
  \selectfont}%
\fi\endgroup%
{\renewcommand{\dashlinestretch}{30}
\begin{picture}(624,339)(0,-10)
\path(12,12)(87,87)
\path(312,312)(237,237)
\path(12,312)(312,12)
\path(312,312)(612,12)
\path(612,312)(537,237)
\path(312,12)(387,87)
\end{picture}
}
 }
    \hspace{-4mm}
  \end{array}
$. Its Khovanov complex is the complex $\Omega_1=\Kh(
  \setlength{\unitlength}{0.5\standardunitlength}
  \begin{array}{c}  \hspace{-2.8mm}
    \raisebox{-2pt}{%
\begingroup\makeatletter\ifx\SetFigFont\undefined%
\gdef\SetFigFont#1#2#3#4#5{%
  \reset@font\fontsize{#1}{#2pt}%
  \fontfamily{#3}\fontseries{#4}\fontshape{#5}%
  \selectfont}%
\fi\endgroup%
{\renewcommand{\dashlinestretch}{30}
\begin{picture}(624,339)(0,-10)
\path(12,12)(87,87)
\path(312,312)(237,237)
\path(12,312)(312,12)
\path(312,312)(612,12)
\path(612,312)(537,237)
\path(312,12)(387,87)
\end{picture}
}
 }
    \hspace{-4mm}
  \end{array}
)$
appearing below. Here and later we use the same notational conventions
as in~\cite{Bar-Natan:Cobordisms}. Thus the height zero object in the
complex is underlined, and $
  \setlength{\unitlength}{0.5\standardunitlength}
  \begin{array}{c}  \hspace{-2.8mm}
    \raisebox{-2pt}{%
\begingroup\makeatletter\ifx\SetFigFont\undefined%
\gdef\SetFigFont#1#2#3#4#5{%
  \reset@font\fontsize{#1}{#2pt}%
  \fontfamily{#3}\fontseries{#4}\fontshape{#5}%
  \selectfont}%
\fi\endgroup%
{\renewcommand{\dashlinestretch}{30}
\begin{picture}(316,331)(0,-10)
\put(158.000,420.500){\arc{375.000}{0.6435}{2.4981}}
\put(158.000,-104.500){\arc{375.000}{3.7851}{5.6397}}
\path(158,233)(158,83)
\end{picture}
}
 }
    \hspace{-4mm}
  \end{array}
$ is the saddle cobordism whose domain
is $
  \setlength{\unitlength}{0.5\standardunitlength}
  \begin{array}{c}  \hspace{-2.8mm}
    \raisebox{-2pt}{%
\begingroup\makeatletter\ifx\SetFigFont\undefined%
\gdef\SetFigFont#1#2#3#4#5{%
  \reset@font\fontsize{#1}{#2pt}%
  \fontfamily{#3}\fontseries{#4}\fontshape{#5}%
  \selectfont}%
\fi\endgroup%
{\renewcommand{\dashlinestretch}{30}
\begin{picture}(316,331)(0,-10)
\put(158.000,420.500){\arc{375.000}{0.6435}{2.4981}}
\put(158.000,-104.500){\arc{375.000}{3.7851}{5.6397}}
\end{picture}
}
 }
    \hspace{-4mm}
  \end{array}
$ and whose range is $
  \setlength{\unitlength}{0.5\standardunitlength}
  \begin{array}{c}  \hspace{-2.8mm}
    \raisebox{-2pt}{%
\begingroup\makeatletter\ifx\SetFigFont\undefined%
\gdef\SetFigFont#1#2#3#4#5{%
  \reset@font\fontsize{#1}{#2pt}%
  \fontfamily{#3}\fontseries{#4}\fontshape{#5}%
  \selectfont}%
\fi\endgroup%
{\renewcommand{\dashlinestretch}{30}
\begin{picture}(316,331)(0,-10)
\put(-104.500,158.000){\arc{375.000}{5.3559}{7.2105}}
\put(420.500,158.000){\arc{375.000}{2.2143}{4.0689}}
\end{picture}
}
 }
    \hspace{-4mm}
  \end{array}
$:

\[ \Omega_1:\qquad
  \setlength{\unitlength}{0.75\standardunitlength}
  \begin{array}{c}
    {\begingroup\makeatletter\ifx\SetFigFont\undefined%
\gdef\SetFigFont#1#2#3#4#5{%
  \reset@font\fontsize{#1}{#2pt}%
  \fontfamily{#3}\fontseries{#4}\fontshape{#5}%
  \selectfont}%
\fi\endgroup%
{\renewcommand{\dashlinestretch}{30}
\begin{picture}(7900,1685)(0,-10)
\put(1887.000,1174.500){\arc{375.000}{3.7851}{5.6397}}
\put(2187.000,1699.000){\arc{374.401}{0.6414}{2.5002}}
\put(2187.000,1175.000){\arc{374.401}{3.7830}{5.6418}}
\path(2187,1362)(2187,1512)
\put(1887.000,1099.000){\arc{374.401}{0.6414}{2.5002}}
\put(1887.000,575.000){\arc{374.401}{3.7830}{5.6418}}
\path(1887,762)(1887,912)
\put(2187.000,1099.500){\arc{375.000}{0.6435}{2.4981}}
\put(2187.000,574.500){\arc{375.000}{3.7851}{5.6397}}
\put(3462.000,1137.000){\arc{3750.000}{2.8578}{3.4254}}
\put(612.000,1137.000){\arc{3750.000}{5.9994}{6.5670}}
\put(4587.000,1099.000){\arc{374.401}{0.6414}{2.5002}}
\put(4587.000,575.000){\arc{374.401}{3.7830}{5.6418}}
\path(4587,762)(4587,912)
\put(4624.500,837.000){\arc{375.000}{5.3559}{7.2105}}
\put(5149.500,837.000){\arc{375.000}{2.2143}{4.0689}}
\put(5374.500,837.000){\arc{375.000}{5.3559}{7.2105}}
\put(5899.500,837.000){\arc{375.000}{2.2143}{4.0689}}
\put(5937.000,1099.000){\arc{374.401}{0.6414}{2.5002}}
\put(5937.000,575.000){\arc{374.401}{3.7830}{5.6418}}
\path(5937,762)(5937,912)
\put(4662.000,837.000){\arc{750.000}{2.4981}{3.7851}}
\put(5862.000,837.000){\arc{750.000}{5.6397}{6.9267}}
\path(5487,837)(5337,837)
\put(537.000,799.500){\arc{375.000}{0.6435}{2.4981}}
\put(1887.000,1699.500){\arc{375.000}{0.6435}{2.4981}}
\put(537.000,274.500){\arc{375.000}{3.7851}{5.6397}}
\put(1062,462){\makebox(0,0)[b]{\smash{{\mathmode{\{-4\}}}}}}
\put(237.000,799.500){\arc{375.000}{0.6435}{2.4981}}
\put(237.000,274.500){\arc{375.000}{3.7851}{5.6397}}
\path(87,762)(12,762)(12,312)(87,312)
\path(687,762)(762,762)(762,312)(687,312)
\put(3087.000,1099.500){\arc{375.000}{0.6435}{2.4981}}
\put(3087.000,574.500){\arc{375.000}{3.7851}{5.6397}}
\put(3124.500,837.000){\arc{375.000}{5.3559}{7.2105}}
\put(3649.500,837.000){\arc{375.000}{2.2143}{4.0689}}
\put(3387.000,499.500){\arc{375.000}{0.6435}{2.4981}}
\put(3387.000,-25.500){\arc{375.000}{3.7851}{5.6397}}
\put(2824.500,237.000){\arc{375.000}{5.3559}{7.2105}}
\put(3349.500,237.000){\arc{375.000}{2.2143}{4.0689}}
\put(6499.500,537.000){\arc{375.000}{5.3559}{7.2105}}
\put(7024.500,537.000){\arc{375.000}{2.2143}{4.0689}}
\put(6799.500,537.000){\arc{375.000}{5.3559}{7.2105}}
\put(7324.500,537.000){\arc{375.000}{2.2143}{4.0689}}
\path(1362,537)(2712,537)
\path(2592.000,507.000)(2712.000,537.000)(2592.000,567.000)
\path(4212,537)(6312,537)
\path(6192.000,507.000)(6312.000,537.000)(6192.000,567.000)
\path(2937,1062)(2862,1062)(2862,12)(2937,12)
\path(3537,1062)(3612,1062)(3612,12)(3537,12)
\path(6537,762)(6462,762)(6462,312)(6537,312)
\path(7287,762)(7362,762)(7362,312)(7287,312)
\path(6462,237)(7887,237)
\put(3912,462){\makebox(0,0)[b]{\smash{{\mathmode{\{-3\}}}}}}
\put(7662,462){\makebox(0,0)[b]{\smash{{\mathmode{\{-2\}}}}}}
\end{picture}
} } 
  \end{array}
 \]

Only one of the objects in $\Omega_1$ contains a loop --- the last
one. Delooping it using Lemma~\ref{lem:Delooping} we get the complex
$\Omega_2$, which is isomorphic to $\Omega_1$ (we've also replaced some of
the smoothings and cobordisms appearing in $\Omega_1$ with isotopic ones):

\[ \Omega_2:\qquad
  \setlength{\unitlength}{0.75\standardunitlength}
  \begin{array}{c}
    {\begingroup\makeatletter\ifx\SetFigFont\undefined%
\gdef\SetFigFont#1#2#3#4#5{%
  \reset@font\fontsize{#1}{#2pt}%
  \fontfamily{#3}\fontseries{#4}\fontshape{#5}%
  \selectfont}%
\fi\endgroup%
{\renewcommand{\dashlinestretch}{30}
\begin{picture}(6924,1760)(0,-10)
\put(6349.500,912.000){\arc{375.000}{2.2143}{4.0689}}
\put(5824.500,312.000){\arc{375.000}{5.3559}{7.2105}}
\put(6349.500,312.000){\arc{375.000}{2.2143}{4.0689}}
\put(4024.500,1512.000){\arc{375.000}{5.3559}{7.2105}}
\put(4549.500,1512.000){\arc{375.000}{2.2143}{4.0689}}
\put(4024.500,912.000){\arc{375.000}{5.3559}{7.2105}}
\put(4549.500,912.000){\arc{375.000}{2.2143}{4.0689}}
\put(4212,912){\blacken\ellipse{60}{60}}
\put(4212,912){\ellipse{60}{60}}
\put(4774.500,1512.000){\arc{375.000}{5.3559}{7.2105}}
\put(5299.500,1512.000){\arc{375.000}{2.2143}{4.0689}}
\put(4774.500,912.000){\arc{375.000}{5.3559}{7.2105}}
\put(5299.500,912.000){\arc{375.000}{2.2143}{4.0689}}
\put(5112,912){\blacken\ellipse{60}{60}}
\put(5112,912){\ellipse{60}{60}}
\put(5862.000,1212.000){\arc{3750.000}{2.8578}{3.4254}}
\put(3462.000,1212.000){\arc{3750.000}{5.9994}{6.5670}}
\path(4887,912)(4737,912)
\path(4887,1512)(4737,1512)
\put(2524.500,912.000){\arc{375.000}{5.3559}{7.2105}}
\put(3049.500,912.000){\arc{375.000}{2.2143}{4.0689}}
\put(2524.500,312.000){\arc{375.000}{5.3559}{7.2105}}
\put(3049.500,312.000){\arc{375.000}{2.2143}{4.0689}}
\put(237.000,874.500){\arc{375.000}{0.6435}{2.4981}}
\put(5824.500,912.000){\arc{375.000}{5.3559}{7.2105}}
\put(237.000,349.500){\arc{375.000}{3.7851}{5.6397}}
\put(4662,312){\makebox(0,0)[b]{\smash{{\mathmode{d^{-1}}}}}}
\put(1737.000,1774.000){\arc{374.401}{0.6414}{2.5002}}
\put(1737.000,1250.000){\arc{374.401}{3.7830}{5.6418}}
\path(1737,1437)(1737,1587)
\put(1737.000,1174.000){\arc{374.401}{0.6414}{2.5002}}
\put(1737.000,650.000){\arc{374.401}{3.7830}{5.6418}}
\path(1737,837)(1737,987)
\put(3312.000,1212.000){\arc{3750.000}{2.8578}{3.4254}}
\put(162.000,1212.000){\arc{3750.000}{5.9994}{6.5670}}
\path(1062,612)(2412,612)
\path(2292.000,582.000)(2412.000,612.000)(2292.000,642.000)
\path(387,837)(462,837)(462,387)(387,387)
\path(87,837)(12,837)(12,387)(87,387)
\put(762,537){\makebox(0,0)[b]{\smash{{\mathmode{\{-4\}}}}}}
\put(1737,312){\makebox(0,0)[b]{\smash{{\mathmode{d^{-2}}}}}}
\path(3612,612)(5712,612)
\path(5592.000,582.000)(5712.000,612.000)(5592.000,642.000)
\path(2937,1137)(3012,1137)(3012,87)(2937,87)
\path(5862,12)(6912,12)
\path(5937,1137)(5862,1137)(5862,87)(5937,87)
\path(6837,1137)(6912,1137)(6912,87)(6837,87)
\path(2637,1137)(2562,1137)(2562,87)(2637,87)
\put(3312,537){\makebox(0,0)[b]{\smash{{\mathmode{\{-3\}}}}}}
\put(6612,837){\makebox(0,0)[b]{\smash{{\mathmode{\{-3\}}}}}}
\put(6612,237){\makebox(0,0)[b]{\smash{{\mathmode{\{-1\}}}}}}
\end{picture}
} } 
  \end{array}
 \]

The symbol $
  \setlength{\unitlength}{0.5\standardunitlength}
  \begin{array}{c}  \hspace{-2.8mm}
    \raisebox{-2pt}{%
\begingroup\makeatletter\ifx\SetFigFont\undefined%
\gdef\SetFigFont#1#2#3#4#5{%
  \reset@font\fontsize{#1}{#2pt}%
  \fontfamily{#3}\fontseries{#4}\fontshape{#5}%
  \selectfont}%
\fi\endgroup%
{\renewcommand{\dashlinestretch}{30}
\begin{picture}(316,331)(0,-10)
\put(-104.500,158.000){\arc{375.000}{5.3559}{7.2105}}
\put(420.500,158.000){\arc{375.000}{2.2143}{4.0689}}
\end{picture}
}
 }
    \hspace{-4mm}
  \end{array}
$ here, when appearing as a cobordism, denotes the
identity automorphism of the smoothing $
  \setlength{\unitlength}{0.5\standardunitlength}
  \begin{array}{c}  \hspace{-2.8mm}
    \raisebox{-2pt}{%
\begingroup\makeatletter\ifx\SetFigFont\undefined%
\gdef\SetFigFont#1#2#3#4#5{%
  \reset@font\fontsize{#1}{#2pt}%
  \fontfamily{#3}\fontseries{#4}\fontshape{#5}%
  \selectfont}%
\fi\endgroup%
{\renewcommand{\dashlinestretch}{30}
\begin{picture}(316,331)(0,-10)
\put(-104.500,158.000){\arc{375.000}{5.3559}{7.2105}}
\put(420.500,158.000){\arc{375.000}{2.2143}{4.0689}}
\end{picture}
}
 }
    \hspace{-4mm}
  \end{array}
$. Likewise $
  \setlength{\unitlength}{0.5\standardunitlength}
  \begin{array}{c}  \hspace{-2.8mm}
    \raisebox{-2pt}{%
\begingroup\makeatletter\ifx\SetFigFont\undefined%
\gdef\SetFigFont#1#2#3#4#5{%
  \reset@font\fontsize{#1}{#2pt}%
  \fontfamily{#3}\fontseries{#4}\fontshape{#5}%
  \selectfont}%
\fi\endgroup%
{\renewcommand{\dashlinestretch}{30}
\begin{picture}(316,331)(0,-10)
\put(-104.500,158.000){\arc{375.000}{5.3559}{7.2105}}
\put(420.500,158.000){\arc{375.000}{2.2143}{4.0689}}
\put(83,158){\blacken\ellipse{60}{60}}
\put(83,158){\ellipse{60}{60}}
\end{picture}
}
 }
    \hspace{-4mm}
  \end{array}
$ 
and $
  \setlength{\unitlength}{0.5\standardunitlength}
  \begin{array}{c}  \hspace{-2.8mm}
    \raisebox{-2pt}{%
\begingroup\makeatletter\ifx\SetFigFont\undefined%
\gdef\SetFigFont#1#2#3#4#5{%
  \reset@font\fontsize{#1}{#2pt}%
  \fontfamily{#3}\fontseries{#4}\fontshape{#5}%
  \selectfont}%
\fi\endgroup%
{\renewcommand{\dashlinestretch}{30}
\begin{picture}(316,331)(0,-10)
\put(-104.500,158.000){\arc{375.000}{5.3559}{7.2105}}
\put(420.500,158.000){\arc{375.000}{2.2143}{4.0689}}
\put(233,158){\blacken\ellipse{60}{60}}
\put(233,158){\ellipse{60}{60}}
\end{picture}
}
 }
    \hspace{-4mm}
  \end{array}
$ denote that same cobordism, with an extra dot on the left
(or right) ``curtain''.

The upper left entry in the differential $d^{-1}$ in $\Omega_2$ is
an isomorphism (so is the upper right entry, but one is enough). So
we can apply the first part of Lemma~\ref{lem:GaussianElimination}
taking $\phi$ to be that upper left entry.  The result is the
complex $\Omega_3$ below, which is isomorphic to $\Omega_2$ and
$\Omega_1$. Note that in this case we have no ``$\begin{pmatrix}
\mu & \nu \end{pmatrix}$'' term, and that in this case,
$-\gamma\phi^{-1}\delta=-
  \setlength{\unitlength}{0.5\standardunitlength}
  \begin{array}{c}  \hspace{-2.8mm}
    \raisebox{-2pt}{%
\begingroup\makeatletter\ifx\SetFigFont\undefined%
\gdef\SetFigFont#1#2#3#4#5{%
  \reset@font\fontsize{#1}{#2pt}%
  \fontfamily{#3}\fontseries{#4}\fontshape{#5}%
  \selectfont}%
\fi\endgroup%
{\renewcommand{\dashlinestretch}{30}
\begin{picture}(316,331)(0,-10)
\put(-104.500,158.000){\arc{375.000}{5.3559}{7.2105}}
\put(420.500,158.000){\arc{375.000}{2.2143}{4.0689}}
\put(83,158){\blacken\ellipse{60}{60}}
\put(83,158){\ellipse{60}{60}}
\end{picture}
}
 }
    \hspace{-4mm}
  \end{array}
\circ
  \setlength{\unitlength}{0.5\standardunitlength}
  \begin{array}{c}  \hspace{-2.8mm}
    \raisebox{-2pt}{%
\begingroup\makeatletter\ifx\SetFigFont\undefined%
\gdef\SetFigFont#1#2#3#4#5{%
  \reset@font\fontsize{#1}{#2pt}%
  \fontfamily{#3}\fontseries{#4}\fontshape{#5}%
  \selectfont}%
\fi\endgroup%
{\renewcommand{\dashlinestretch}{30}
\begin{picture}(316,331)(0,-10)
\put(-104.500,158.000){\arc{375.000}{5.3559}{7.2105}}
\put(420.500,158.000){\arc{375.000}{2.2143}{4.0689}}
\end{picture}
}
 }
    \hspace{-4mm}
  \end{array}
\circ(-
  \setlength{\unitlength}{0.5\standardunitlength}
  \begin{array}{c}  \hspace{-2.8mm}
    \raisebox{-2pt}{%
\begingroup\makeatletter\ifx\SetFigFont\undefined%
\gdef\SetFigFont#1#2#3#4#5{%
  \reset@font\fontsize{#1}{#2pt}%
  \fontfamily{#3}\fontseries{#4}\fontshape{#5}%
  \selectfont}%
\fi\endgroup%
{\renewcommand{\dashlinestretch}{30}
\begin{picture}(316,331)(0,-10)
\put(-104.500,158.000){\arc{375.000}{5.3559}{7.2105}}
\put(420.500,158.000){\arc{375.000}{2.2143}{4.0689}}
\end{picture}
}
 }
    \hspace{-4mm}
  \end{array}
)=+
  \setlength{\unitlength}{0.5\standardunitlength}
  \begin{array}{c}  \hspace{-2.8mm}
    \raisebox{-2pt}{%
\begingroup\makeatletter\ifx\SetFigFont\undefined%
\gdef\SetFigFont#1#2#3#4#5{%
  \reset@font\fontsize{#1}{#2pt}%
  \fontfamily{#3}\fontseries{#4}\fontshape{#5}%
  \selectfont}%
\fi\endgroup%
{\renewcommand{\dashlinestretch}{30}
\begin{picture}(316,331)(0,-10)
\put(-104.500,158.000){\arc{375.000}{5.3559}{7.2105}}
\put(420.500,158.000){\arc{375.000}{2.2143}{4.0689}}
\put(83,158){\blacken\ellipse{60}{60}}
\put(83,158){\ellipse{60}{60}}
\end{picture}
}
 }
    \hspace{-4mm}
  \end{array}
$:

\[ \Omega_3:\qquad
  \setlength{\unitlength}{0.75\standardunitlength}
  \begin{array}{c}
    {\begingroup\makeatletter\ifx\SetFigFont\undefined%
\gdef\SetFigFont#1#2#3#4#5{%
  \reset@font\fontsize{#1}{#2pt}%
  \fontfamily{#3}\fontseries{#4}\fontshape{#5}%
  \selectfont}%
\fi\endgroup%
{\renewcommand{\dashlinestretch}{30}
\begin{picture}(6924,1760)(0,-10)
\put(6349.500,912.000){\arc{375.000}{2.2143}{4.0689}}
\put(5824.500,312.000){\arc{375.000}{5.3559}{7.2105}}
\put(6349.500,312.000){\arc{375.000}{2.2143}{4.0689}}
\put(3724.500,1512.000){\arc{375.000}{5.3559}{7.2105}}
\put(4249.500,1512.000){\arc{375.000}{2.2143}{4.0689}}
\put(4399.500,912.000){\arc{375.000}{5.3559}{7.2105}}
\put(4924.500,912.000){\arc{375.000}{2.2143}{4.0689}}
\put(4737,912){\blacken\ellipse{60}{60}}
\put(4737,912){\ellipse{60}{60}}
\put(4999.500,912.000){\arc{375.000}{5.3559}{7.2105}}
\put(5524.500,912.000){\arc{375.000}{2.2143}{4.0689}}
\put(5187,912){\blacken\ellipse{60}{60}}
\put(5187,912){\ellipse{60}{60}}
\put(2524.500,912.000){\arc{375.000}{5.3559}{7.2105}}
\put(3049.500,912.000){\arc{375.000}{2.2143}{4.0689}}
\put(2524.500,312.000){\arc{375.000}{5.3559}{7.2105}}
\put(3049.500,312.000){\arc{375.000}{2.2143}{4.0689}}
\put(237.000,874.500){\arc{375.000}{0.6435}{2.4981}}
\put(237.000,349.500){\arc{375.000}{3.7851}{5.6397}}
\put(1737.000,1174.000){\arc{374.401}{0.6414}{2.5002}}
\put(1737.000,650.000){\arc{374.401}{3.7830}{5.6418}}
\path(1737,837)(1737,987)
\put(3312.000,1212.000){\arc{3750.000}{2.8578}{3.4254}}
\put(162.000,1212.000){\arc{3750.000}{5.9994}{6.5670}}
\put(5824.500,912.000){\arc{375.000}{5.3559}{7.2105}}
\path(1062,612)(2412,612)
\path(2292.000,582.000)(2412.000,612.000)(2292.000,642.000)
\put(3987,837){\makebox(0,0)[b]{\smash{{\mathmode{0}}}}}
\path(387,837)(462,837)(462,387)(387,387)
\path(87,837)(12,837)(12,387)(87,387)
\put(762,537){\makebox(0,0)[b]{\smash{{\mathmode{\{-4\}}}}}}
\put(1737,312){\makebox(0,0)[b]{\smash{{\mathmode{d^{-2}}}}}}
\put(1737,1437){\makebox(0,0)[b]{\smash{{\mathmode{0}}}}}
\put(5562.000,1212.000){\arc{3750.000}{2.8578}{3.4254}}
\put(3762.000,1212.000){\arc{3750.000}{5.9994}{6.5670}}
\path(3612,612)(5712,612)
\path(5592.000,582.000)(5712.000,612.000)(5592.000,642.000)
\path(2937,1137)(3012,1137)(3012,87)(2937,87)
\path(5862,12)(6912,12)
\path(5937,1137)(5862,1137)(5862,87)(5937,87)
\path(6837,1137)(6912,1137)(6912,87)(6837,87)
\path(4512,912)(4362,912)
\path(5037,912)(4887,912)
\path(4962,987)(4962,837)
\path(2637,1137)(2562,1137)(2562,87)(2637,87)
\put(3312,537){\makebox(0,0)[b]{\smash{{\mathmode{\{-3\}}}}}}
\put(6612,837){\makebox(0,0)[b]{\smash{{\mathmode{\{-3\}}}}}}
\put(6612,237){\makebox(0,0)[b]{\smash{{\mathmode{\{-1\}}}}}}
\put(4662,312){\makebox(0,0)[b]{\smash{{\mathmode{d^{-1}}}}}}
\put(4962,1437){\makebox(0,0)[b]{\smash{{\mathmode{0}}}}}
\end{picture}
} } 
  \end{array}
 \]

Dropping the contractible summand as in the second part of
Lemma~\ref{lem:GaussianElimination} we get the complex $\Omega_4$:

\[ \Omega_4:\qquad
  \setlength{\unitlength}{0.75\standardunitlength}
  \begin{array}{c}
    {\begingroup\makeatletter\ifx\SetFigFont\undefined%
\gdef\SetFigFont#1#2#3#4#5{%
  \reset@font\fontsize{#1}{#2pt}%
  \fontfamily{#3}\fontseries{#4}\fontshape{#5}%
  \selectfont}%
\fi\endgroup%
{\renewcommand{\dashlinestretch}{30}
\begin{picture}(6775,860)(0,-10)
\put(237.000,49.500){\arc{375.000}{3.7851}{5.6397}}
\put(1737.000,874.000){\arc{374.401}{0.6414}{2.5002}}
\put(1737.000,350.000){\arc{374.401}{3.7830}{5.6418}}
\path(1737,537)(1737,687)
\put(4099.500,612.000){\arc{375.000}{5.3559}{7.2105}}
\put(4624.500,612.000){\arc{375.000}{2.2143}{4.0689}}
\put(4437,612){\blacken\ellipse{60}{60}}
\put(4437,612){\ellipse{60}{60}}
\put(4699.500,612.000){\arc{375.000}{5.3559}{7.2105}}
\put(5224.500,612.000){\arc{375.000}{2.2143}{4.0689}}
\put(4887,612){\blacken\ellipse{60}{60}}
\put(4887,612){\ellipse{60}{60}}
\path(4212,612)(4062,612)
\path(4737,612)(4587,612)
\path(4662,687)(4662,537)
\put(2524.500,312.000){\arc{375.000}{5.3559}{7.2105}}
\put(3049.500,312.000){\arc{375.000}{2.2143}{4.0689}}
\path(2637,537)(2562,537)(2562,87)(2637,87)
\put(237.000,574.500){\arc{375.000}{0.6435}{2.4981}}
\path(2937,537)(3012,537)(3012,87)(2937,87)
\put(6537,237){\makebox(0,0)[b]{\smash{{\mathmode{\{-1\}}}}}}
\put(3312,237){\makebox(0,0)[b]{\smash{{\mathmode{\{-3\}}}}}}
\put(5749.500,312.000){\arc{375.000}{5.3559}{7.2105}}
\put(6274.500,312.000){\arc{375.000}{2.2143}{4.0689}}
\put(1812.000,612.000){\arc{750.000}{2.4981}{3.7851}}
\put(1662.000,612.000){\arc{750.000}{5.6397}{6.9267}}
\put(4287.000,612.000){\arc{750.000}{2.4981}{3.7851}}
\put(4887.000,612.000){\arc{750.000}{5.6397}{6.9267}}
\path(1062,312)(2412,312)
\path(2292.000,282.000)(2412.000,312.000)(2292.000,342.000)
\path(387,537)(462,537)(462,87)(387,87)
\path(87,537)(12,537)(12,87)(87,87)
\path(3612,312)(5712,312)
\path(5592.000,282.000)(5712.000,312.000)(5592.000,342.000)
\path(5862,537)(5787,537)(5787,87)(5862,87)
\path(6162,537)(6237,537)(6237,87)(6162,87)
\path(5787,12)(6762,12)
\put(762,237){\makebox(0,0)[b]{\smash{{\mathmode{\{-4\}}}}}}
\end{picture}
} } 
  \end{array}
 \]

The complex $\Omega_4$ contains fewer objects than $\Omega_1$, hence
from a computational perspective it is indeed simpler. The savings, about
25\%, may not appear to be much, but remember ---
\begin{itemize}
\item This 25\% will be compounded with whatever savings we will incur
later in the computation.
\item This is just a sample. In real life the savings are greater.
\end{itemize}

In a similar manner we can compute and simplify the Khovanov complex of the
right half of the Figure Eight knot, the tangle $T_R$. The result is the
complex $\Psi_4$:

\[ \Psi_4:\qquad
  \setlength{\unitlength}{0.75\standardunitlength}
  \begin{array}{c}
    {\begingroup\makeatletter\ifx\SetFigFont\undefined%
\gdef\SetFigFont#1#2#3#4#5{%
  \reset@font\fontsize{#1}{#2pt}%
  \fontfamily{#3}\fontseries{#4}\fontshape{#5}%
  \selectfont}%
\fi\endgroup%
{\renewcommand{\dashlinestretch}{30}
\begin{picture}(6796,860)(0,-10)
\put(6349.500,312.000){\arc{375.000}{2.2143}{4.0689}}
\put(1812.000,350.000){\arc{374.401}{3.7830}{5.6418}}
\put(1812.000,874.000){\arc{374.401}{0.6414}{2.5002}}
\put(1812,687){\blacken\ellipse{60}{60}}
\put(1812,687){\ellipse{60}{60}}
\put(2412.000,350.000){\arc{374.401}{3.7830}{5.6418}}
\put(2412.000,874.000){\arc{374.401}{0.6414}{2.5002}}
\put(2412,537){\blacken\ellipse{60}{60}}
\put(2412,537){\ellipse{60}{60}}
\put(237.000,574.500){\arc{375.000}{0.6435}{2.4981}}
\put(237.000,49.500){\arc{375.000}{3.7851}{5.6397}}
\put(3537.000,574.500){\arc{375.000}{0.6435}{2.4981}}
\put(3537.000,49.500){\arc{375.000}{3.7851}{5.6397}}
\put(5037.000,350.000){\arc{374.401}{3.7830}{5.6418}}
\put(5037.000,874.000){\arc{374.401}{0.6414}{2.5002}}
\path(5037,687)(5037,537)
\put(4962.000,612.000){\arc{750.000}{5.6397}{6.9267}}
\put(1737.000,612.000){\arc{750.000}{2.4981}{3.7851}}
\put(5824.500,312.000){\arc{375.000}{5.3559}{7.2105}}
\put(2337.000,612.000){\arc{750.000}{5.6397}{6.9267}}
\put(4062,237){\makebox(0,0)[b]{\smash{{\mathmode{\{3\}}}}}}
\put(5112.000,612.000){\arc{750.000}{2.4981}{3.7851}}
\path(5937,537)(5862,537)(5862,87)(5937,87)
\path(6237,537)(6312,537)(6312,87)(6237,87)
\path(4362,312)(5712,312)
\path(5592.000,282.000)(5712.000,312.000)(5592.000,342.000)
\path(1062,312)(3162,312)
\path(3042.000,282.000)(3162.000,312.000)(3042.000,342.000)
\path(1662,612)(1512,612)
\path(2187,612)(2037,612)
\path(2112,687)(2112,537)
\path(12,12)(987,12)
\path(387,537)(462,537)(462,87)(387,87)
\path(87,537)(12,537)(12,87)(87,87)
\path(3687,537)(3762,537)(3762,87)(3687,87)
\path(3387,537)(3312,537)(3312,87)(3387,87)
\put(6612,237){\makebox(0,0)[b]{\smash{{\mathmode{\{4\}}}}}}
\put(762,237){\makebox(0,0)[b]{\smash{{\mathmode{\{1\}}}}}}
\end{picture}
} } 
  \end{array}
 \]

Next we have to take the ``tensor product'' of $\Omega_4$ with $\Psi_4$
using the same side-by-side composition used to make the Figure Eight knot
out of $T_L$ and $T_R$. The result is the double complex $\Upsilon_1$ below.
To prevent clutter we've done away with the vector and matrix brackets;
also, notice the signs in the middle row, sprinkled as commonly practiced
on tensor products.

\[ 
  \setlength{\unitlength}{0.75\standardunitlength}
  \begin{array}{c}
    {\input figs/Upsilon1.eepic } 
  \end{array}
 \]

$\Upsilon_1$ may seem monstrous, but yet, it contains just 9 objects, as
compared to the 16 objects one sees in a direct computation of the Khovanov
complex of a 4-crossing knot. Anyway, it is best to re-write $\Upsilon_1$ a
bit before proceeding. The smoothings and cobordisms can be smoothed out,
and dots can be moved around cobordisms so as to cancel the four
differences on the lower left of the $\Upsilon_1$ diagram. The result is
$\Upsilon_2$ below.

\[ \Upsilon_2:\qquad
  \setlength{\unitlength}{0.75\standardunitlength}
  \begin{array}{c}
    {\input figs/Upsilon2.eepic } 
  \end{array}
 \]

If we were a computer program we would have now ``flattened'' the double
complex $\Upsilon_2$ to a single complex of the schematic form
$(\cdot)\to(:)\to(\vdots)\to(:)\to(\cdot)$. Such ``single'' complexes are
more easily manipulated on a computer. But it is unlikely this paper will
ever be appreciated by anything but humans. So we'll stick to the more
readable double complex form while remembering that we really have just one
differential going south and east.

Anyway, the next step is to replace every loop in every object in
$\Upsilon_2$ with a pair of (degree-shifted) empty sets, as in
Lemma~\ref{lem:Delooping}, while replacing the differentials with their
compositions with the explicit isomorphisms of~\eqref{eq:Delooping}. But
there's nothing but loops in $\Upsilon_2$, so we are left with a complex
$\Upsilon_3$ in which all the objects are degree-shifted empty sets and all
the morphisms are (matrices of) scalar multiples of the empty cobordism
(note that modulo the relations in~\eqref{eq:LocalRelations} all closed
surfaces reduce to scalars):

\[
  \def\s{\scriptstyle}
  \Upsilon_3:\qquad\xymatrix@C=25mm@R=22mm{
    \left[{\begin{matrix}
      \emptyset\{-5\} \\ \emptyset\{-3\} \\ \emptyset\{-3\} \\ \emptyset\{-1\}
    \end{matrix}}\right]
    \ar[r]^{\s\left({\begin{matrix}
      \s 0 &\s 1 &\s -1 &\s 0 \\\s 0 &\s 0 &\s 0 &\s -1 \\
      \s 0 &\s 0 &\s 0 &\s 1 \\\s 0 &\s 0 &\s 0 &\s 0 
    \end{matrix}}\right)}
    \ar[d]^{\s\left({\begin{matrix}
      \s 0 &\s 1 &\s 1 &\s 0 \\\s 0 &\s 0 &\s 0 &\s 1
    \end{matrix}}\right)}_m
  &
    \left[{\begin{matrix}
      \emptyset\{-3\} \\ \emptyset\{-1\} \\ \emptyset\{-1\} \\ \emptyset\{1\}
    \end{matrix}}\right]
    \ar[r]^{\s\left({\begin{matrix}
      \s 0 &\s 1 &\s 1 &\s 0 \\\s 0 &\s 0 &\s 0 &\s 1
    \end{matrix}}\right)}_m
    \ar[d]^(0.4){\s\left({\begin{matrix}
      \s 0 &\s 1 &\s 1 &\s 0 \\\s 0 &\s 0 &\s 0 &\s 1
    \end{matrix}}\right)}_(0.4)m
  &
    \underline{\left[{\begin{matrix}
      \emptyset\{-1\} \\ \emptyset\{1\}
    \end{matrix}}\right]}
    \ar[d]^{\s\left({\begin{matrix}
      \s 1 &\s 0 \\\s 0 &\s 1 \\\s 0 &\s 1 \\\s 0 &\s 0
    \end{matrix}}\right)}_\Delta
  \\
    \left[{\begin{matrix}
      \emptyset\{-3\} \\ \emptyset\{-1\}
    \end{matrix}}\right]
    \ar[r]^0 \ar[d]^0
  &
    \underline{\left[{\begin{matrix}
      \emptyset\{-1\} \\ \emptyset\{1\}
    \end{matrix}}\right]}
    \ar[r]^{\s\left({\begin{matrix}
      \s -1 &\s 0 \\\s 0 &\s -1 \\\s 0 &\s -1 \\\s 0 &\s 0
    \end{matrix}}\right)}_{-\Delta}
    \ar[d]^0
  &
    \left[{\begin{matrix}
      \emptyset\{-1\} \\ \emptyset\{1\} \\ \emptyset\{1\} \\ \emptyset\{3\}
    \end{matrix}}\right]
    \ar[d]^{\s\left({\begin{matrix}
      \s 0 &\s 1 &\s -1 &\s 0 \\\s 0 &\s 0 &\s 0 &\s -1 \\
      \s 0 &\s 0 &\s 0 &\s 1 \\\s 0 &\s 0 &\s 0 &\s 0
    \end{matrix}}\right)}
  \\
    \underline{\left[{\begin{matrix}
      \emptyset\{-1\} \\ \emptyset\{1\}
    \end{matrix}}\right]}
    \ar[r]^0
  & 
    \left[{\begin{matrix}
      \emptyset\{1\} \\ \emptyset\{3\}
    \end{matrix}}\right]
    \ar[r]^{\s\left({\begin{matrix}
      \s 1 &\s 0 \\\s 0 &\s 1 \\\s 0 &\s 1 \\\s 0 &\s 0
    \end{matrix}}\right)}_\Delta
  &
    \left[{\begin{matrix}
      \emptyset\{1\} \\ \emptyset\{3\} \\ \emptyset\{3\} \\ \emptyset\{5\}
    \end{matrix}}\right]
  }
\]

Note in passing that the matrices $m$ and $\Delta$ appearing
here are the matrices representing the product and the
co-product of~\cite{Khovanov:Categorification} relative to
the basis $(X,1)$ used there (or the basis $(v_-, v_+)$ used
in~\cite{Bar-Natan:Categorification}). This is essentially the content
of~\cite[Section~9.1]{Bar-Natan:Cobordisms}.

\parpic[r]{$\displaystyle
  \xymatrix{
    \left[\emptyset\{-5\}\right] \ar[r]^0 \ar[d]^0
  &
    \left[\emptyset\{-1\}\right] \ar[r]^0 \ar[d]^0
  &
    \underline{\left[\right]} \ar[d]^0
  \\
    \left[\right] \ar[r]^0 \ar[d]^0
  &
    \underline{\left[\right]} \ar[r]^0 \ar[d]^0
  &
    \left[\emptyset\{1\}\right] \ar[d]^0
  \\
    \underline{\left[{\begin{matrix}
      \emptyset\{-1\} \\ \emptyset\{1\}
    \end{matrix}}\right]}
    \ar[r]^0
  & 
    \left[\right] \ar[r]^0
  &
    \left[\emptyset\{5\}\right]
  }
$}
We can now apply Lemma~\ref{lem:GaussianElimination} repeatedly to
$\Upsilon_3$, until no invertible entries remain in any of the matrices.
Over $\bbQ$ any non-zero number is invertible so our process stops
when all matrices are $0$. Thus a lengthy iterated application of
Lemma~\ref{lem:GaussianElimination} stops at the complex $\Upsilon_4$
shown on the right (a human could save quite a lot by being clever,
but that's not our point here).

Belatedly flattening the double complex $\Upsilon_4$ we arrive at our final
answer, the complex 
\picskip{0}
\begin{equation} \label{eq:Upsilon5}
  \Upsilon_5:\qquad\xymatrix{
    \left[\emptyset\{-5\}\right] \ar[r]^0 &
    \left[\emptyset\{-1\}\right] \ar[r]^0 &
    \underline{\left[{\begin{matrix}
        \emptyset\{-1\} \\ \emptyset\{1\}
      \end{matrix}}\right]}
      \ar[r]^0 &
    \left[\emptyset\{1\}\right] \ar[r]^0 &
    \left[\emptyset\{5\}\right]
  }.
\end{equation}

To recover homology groups out of $\Upsilon_5$ we need to apply to it
some functor $\calF$ taking $\Cobdl^3$ to graded vector spaces, and
then take homology. The latter step (computing the homology) is the do
nothing operation as all differentials are $0$. Typically (e.g., as
in~\cite[Section 11.2]{Bar-Natan:Cobordisms}) the functor $\calF$ maps
the empty smoothing to the one dimensional vector space (call it
$\bbQ$) at degree $0$. And so we can read directly from~\eqref{eq:Upsilon5}
that the ``conventional'' Khovanov homology over $\bbQ$ of the figure
eight knot is $6$-dimensional with generators at bidegrees $(-2,-5)$,
$(-1,-1)$, $(0,-1)$, $(0,1)$, $(1,1)$ and $(2,5)$.

\section{A faster algorithm} \label{sec:Faster}

%\begin{floatingfigure}[r]{60mm}
%\hspace{-5mm}\includegraphics[width=60mm]{figs/L11n424Scan.ps}
%\vskip -2mm
%\caption{
%  Scanning a knot
%} \label{fig:L11n424Scan}
%\end{floatingfigure}
\parpic[l]{$\includegraphics[width=60mm]{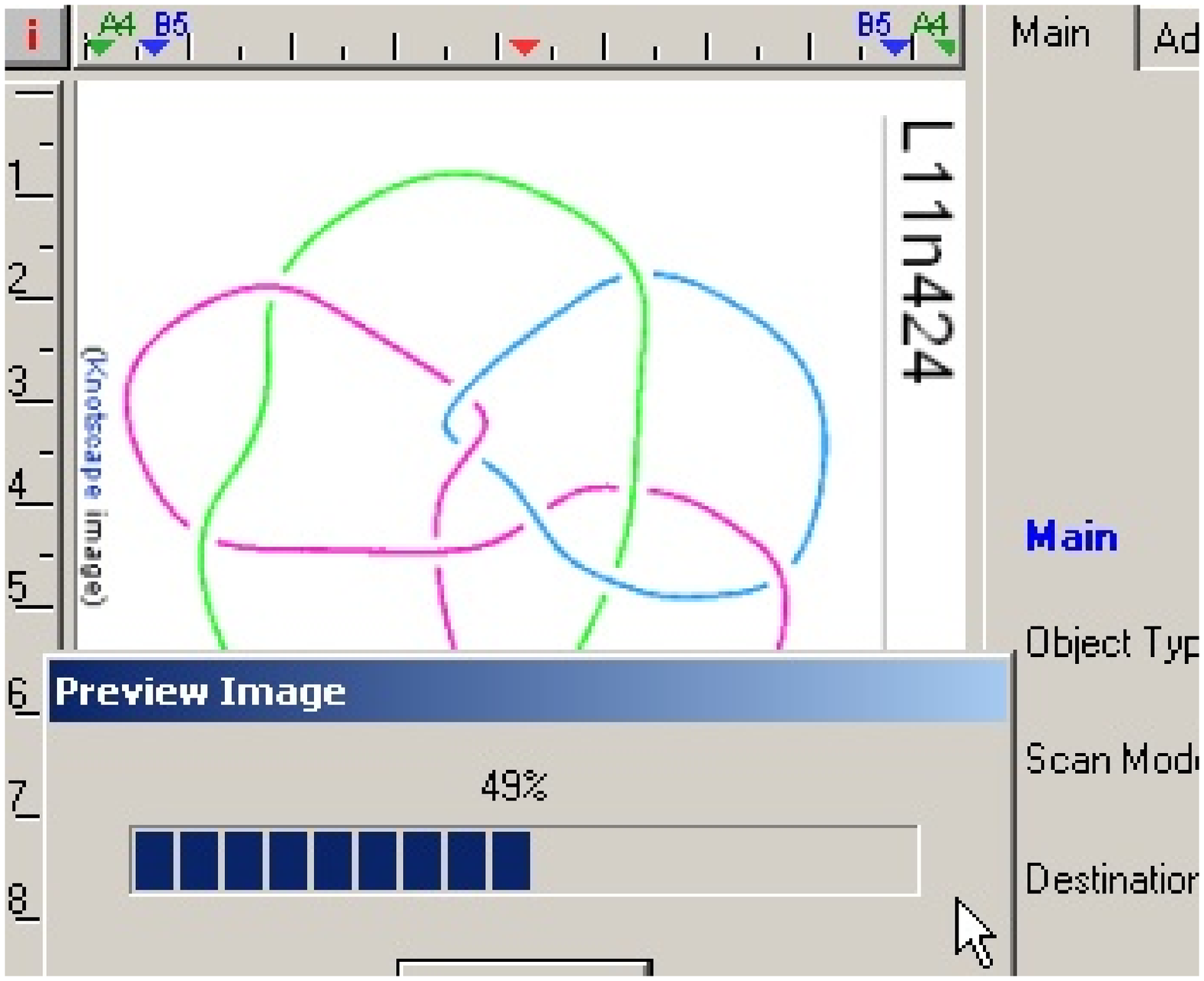}$}
It turns out that there is a somewhat better way to assemble our two main
tools into a running algorithm. Instead of computing each ``half knot'' and
combining the result, ``scan'' the knot from left to right (or top to bottom,
as in the picture on the left), adding one crossing at a time. After
each crossing is added, use delooping (Lemma~\ref{lem:Delooping}) and
Gaussian elimination (Lemma~\ref{lem:GaussianElimination}) to simplify the
result.

Since we have no rigorous estimate, we have to settle with a possibly
naive estimation of the complexity of the ``divide an conquer''
algorithm and of the new ``scanning'' algorithm. In both cases the
bottleneck ought to be where the knot is widest; if the width at the
widest cut is $W$, we expect the complex corresponding to either half
of the knot to be of size $C^W$ for some $C$. In the ``divide an conquer''
algorithm we need to ``tensor multiply'' two such complexes so for a little
while we have to hold a double complex of total size of about
$(C^W)^2=C^{2W}$ (compare with the complex $\Upsilon_1$ above). In the
``scanning'' algorithm we only need to tensor multiply the complex for the
left half-knot with the complex of a single crossing, whose size is $2$. So
we only see a double complex of size $2C^W$ before we get the chance
to simplify again.

\vskip 4mm

\section{Computer programs} \label{sec:ComputerPrograms}

%\begin{floatingfigure}[r]{60mm}
%\hspace{-5mm}\includegraphics[width=60mm]{figs/T87.ps}
%\vskip -2mm
%\caption{
%  The $(8,7)$ torus knot $T_{8,7}$
%} \label{fig:L11n424Scan}
%\end{floatingfigure}
\parpic[r]{$\includegraphics[width=40mm]{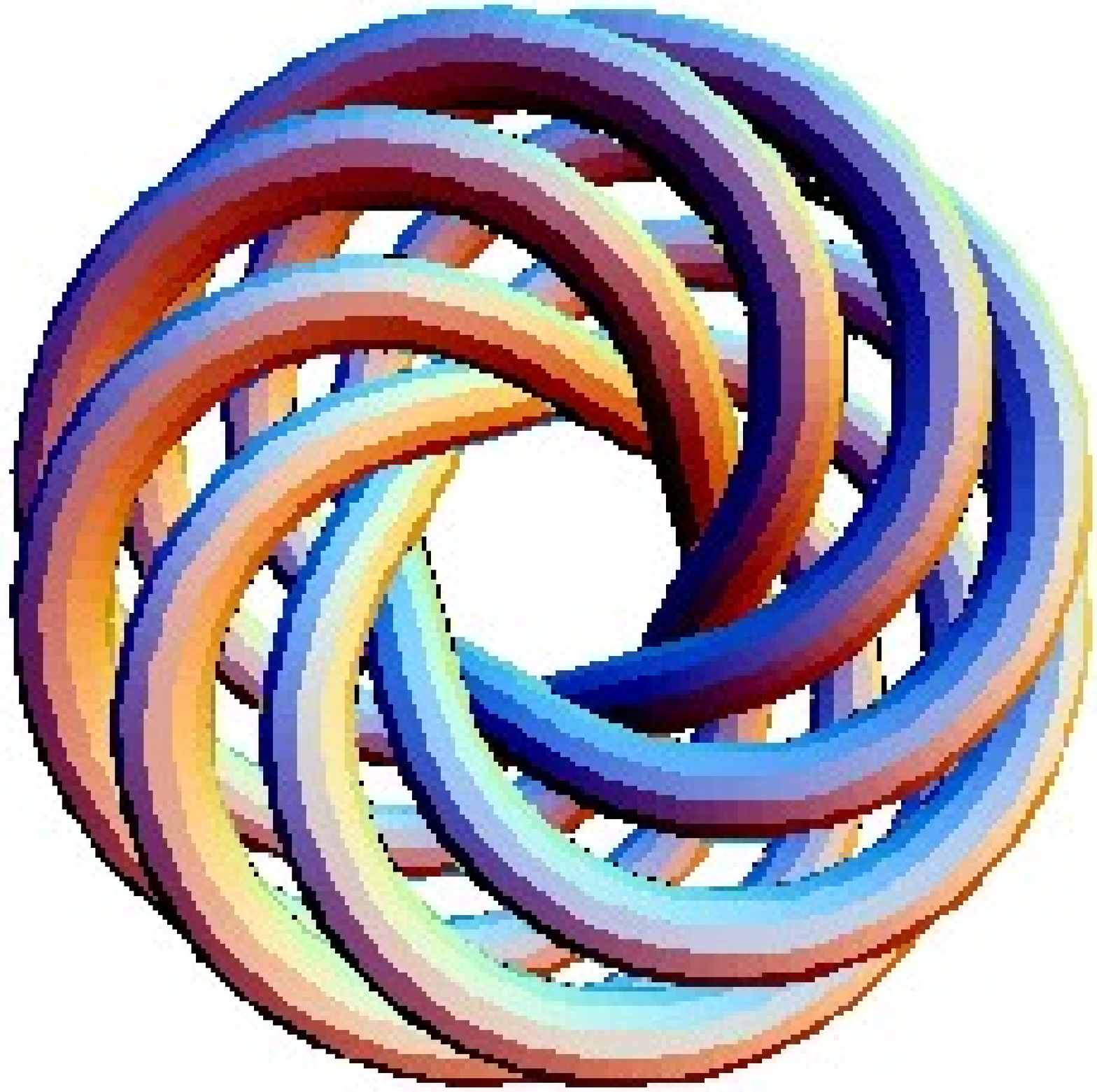}\hspace{-5mm}T_{8,7}$}
Two implementations of the ``scanning'' algorithm are available, both
as a part of the knot theory package {\tt
KnotTheory`}~\cite{KnotTheory}. The first one, {\tt FastKh}, was
written by the author just to test the principle, with no attempt at
optimization. Yet for example, it was able to compute the Khovanov
homology of the $35$-crossing $(7,6)$ torus knot $T_{7,6}$ in about one
day of work; a task that would have taken about a $1,000$ years without
the use of tangles. The second one, {\tt JavaKh} was written by Jeremy
Green, a summer student of the author's at the University of Toronto,
in the summer of 2005. The computation leading to
Table~\ref{table:T87}, for example, of the Khovanov homology of the
$48$-crossing $(8,7)$ torus knot $T_{8,7}$, takes a few minutes using
{\tt JavaKh}. The $(r,j)$ entry of that table contains the degree
$2r+j$ piece of the $r$th integral Khovanov homology
$\calG_{2r+j}H^r(T_{8,7})$ of $T_{8,7}$.

\begin{table}
\def\s{\scriptstyle}
\begin{center} \begin{tabular}{|c||c|c|c|c|c|c|c|c|c|c|c|}
 \hline
 &$\s j=23$ & $\s j=25$ & $\s j=27$ & $\s j=29$ & $\s j=31$ & $\s j=33$ &
$\s \
j=35$ & $\s j=37$ & $\s j=39$ & $\s j=41$ & $\s j=43$ \\
 \hline\hline
  $\s r=0$ & & & & & & & & & & $\s\bbZ$ & $\s\bbZ$ \\
  \hline
  $\s r=1$ & & & & & & & & & & & \\
  \hline
  $\s r=2$ & & & & & & & & & & $\s\bbZ$ & \\
  \hline
  $\s r=3$ & & & & & & & & & & $\s\bbZ_2$ & $\s\bbZ$ \\
  \hline
  $\s r=4$ & & & & & & & & & $\s\bbZ$ & $\s\bbZ$ & \\
  \hline
  $\s r=5$ & & & & & & & & & & $\s\bbZ$ & $\s\bbZ$ \\
  \hline
  $\s r=6$ & & & & & & & & $\s\bbZ$ & $\s\bbZ$ & & \\
  \hline
  $\s r=7$ & & & & & & & & $\s\bbZ_2$ & $\s\bbZ\oplus\bbZ_2$ & $\s\bbZ$ &
\\
  \hline
  $\s r=8$ & & & & & & & $\s\bbZ$ & $\s\bbZ^2$ & & & \\
  \hline
  $\s r=9$ & & & & & & & & $\s\bbZ\oplus\bbZ_2$ & $\s\bbZ^2$ & & \\
  \hline
  $\s r=10$ & & & & & & $\s\bbZ$ & $\s\bbZ^2$ & $\s\bbZ_2$ & $\s\bbZ_2$ & &
\
\\
  \hline
  $\s r=11$ & & & & & & $\s\bbZ_2$ & $\s\bbZ\oplus\bbZ_2^2$ & $\s\bbZ^3$ &
& \
& \\
  \hline
  $\s r=12$ & & & & & $\s\bbZ$ & $\s\bbZ^3$ & $\s\bbZ$ & \
$\s\bbZ_2\oplus\bbZ_5$ & $\s\bbZ$ & & \\
  \hline
  $\s r=13$ & & & & & & $\s\bbZ\oplus\bbZ_2^2$ & $\s\bbZ^4\oplus\bbZ_2$ &
$\s\
\bbZ$ & & & \\
  \hline
  $\s r=14$ & & & & & $\s\bbZ^2$ & $\s\bbZ\oplus\bbZ_2$ & $\s\bbZ_2^2$ & \
$\s\bbZ$ & & & \\
  \hline
  $\s r=15$ & & & & & $\s\bbZ_2^2\oplus\bbZ_7$ & $\s\bbZ^4\oplus\bbZ_2$ &
$\s\
\bbZ^2$ & & & & \\
  \hline
  $\s r=16$ & & & & $\s\bbZ^2$ & $\s\bbZ^2$ & $\s\bbZ_2^2$ & \
$\s\bbZ^2\oplus\bbZ_2$ & $\s\bbZ$ & & & \\
  \hline
  $\s r=17$ & & & & $\s\bbZ_2$ & $\s\bbZ^3\oplus\bbZ_2$ & $\s\bbZ^3$ & & &
& \
& \\
  \hline
  $\s r=18$ & & & $\s\bbZ$ & $\s\bbZ^2\oplus\bbZ_2$ & $\s\bbZ_2^2$ & \
$\s\bbZ\oplus\bbZ_4$ & $\s\bbZ$ & & & & \\
  \hline
  $\s r=19$ & & & $\s\bbZ_2$ & $\s\bbZ^2\oplus\bbZ_2^2$ & $\s\bbZ^3$ & & &
& \
& & \\
  \hline
  $\s r=20$ & & $\s\bbZ$ & $\s\bbZ^2$ & $\s\bbZ_2^2$ &
$\s\bbZ\oplus\bbZ_2^2$ \
& $\s\bbZ$ & & & & & \\
  \hline
  $\s r=21$ & & & $\s\bbZ\oplus\bbZ_4\oplus\bbZ_7$ & $\s\bbZ^3\oplus\bbZ_2$
& \
$\s\bbZ_2$ & & & & & & \\
  \hline
  $\s r=22$ & & $\s\bbZ$ & $\s\bbZ_2$ & \
$\s\bbZ_2\oplus\bbZ_4\oplus\bbZ_5\oplus\bbZ_7$ & $\s\bbZ$ & & & & & & \\
  \hline
  $\s r=23$ & & $\s\bbZ_2\oplus\bbZ_7$ & $\s\bbZ^2$ & $\s\bbZ_2$ &
$\s\bbZ_2$ \
& & & & & & \\
  \hline
  $\s r=24$ & $\s\bbZ$ & $\s\bbZ$ & $\s\bbZ_2\oplus\bbZ_5\oplus\bbZ_7$ & \
$\s\bbZ$ & & & & & & & \\
  \hline
  $\s r=25$ & & $\s\bbZ$ & $\s\bbZ$ & $\s\bbZ_5$ & & & & & & & \\
  \hline
  $\s r=26$ & & & $\s\bbZ_2\oplus\bbZ_3$ & $\s\bbZ_2$ & & & & & & & \\
  \hline
\end{tabular} \end{center}
\caption{The integral Khovanov homology groups $\calG_{2r+j}H^r(T_{8,7})$.}
\label{table:T87}
\end{table}

\section{The Reidemeister moves} \label{sec:Reidemeister}

One further advantage of the ability to compute with tangles is that
the proof of up-to-homotopy invariance of the Khovanov complex under
Reidemeister move becomes mechanical. All that one has to do is to
compute and simplify the complexes corresponding to each side of each
Reidemeister move. One gets the same result for both sides of any given
Reidemeister move and hence invariance is proven\footnote{In fact, the
computer programs discussed above can carry out these computations,
turning them literally mechanical. Though the presently available
front ends for those programs are only configured to take knots as
inputs.}.

\begin{figure}
\[ \includegraphics[width=6.5in]{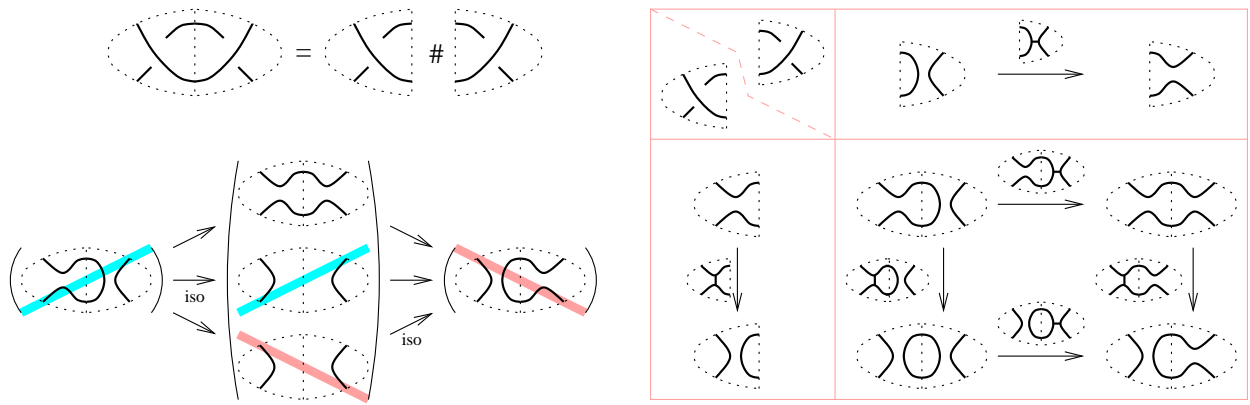} \]
\caption{Proof of invariance under $R2$.}
\label{fig:R2}
\end{figure}

As an example, Figure~\ref{fig:R2} proves invariance under the second
Reidemeister move $R2$ in summary form. Start on the upper left, where
the ``hard side'' of $R2$ is presented as the side-by-side product of
two single-crossing tangles. Move on to the right where the
corresponding double complex is shown, and back to the bottom left
where that complex is flattened, delooped and the two isomorphisms that
appear in the result are canceled out. What remains is a single-entry
complex equal to the ``easy side'' of $R2$.  For simplicity degree
shifts are ignored in Figure~\ref{fig:R2}.

\[ 
  \setlength{\unitlength}{0.5\standardunitlength}
  \begin{array}{c}
    {\begingroup\makeatletter\ifx\SetFigFont\undefined%
\gdef\SetFigFont#1#2#3#4#5{%
  \reset@font\fontsize{#1}{#2pt}%
  \fontfamily{#3}\fontseries{#4}\fontshape{#5}%
  \selectfont}%
\fi\endgroup%
{\renewcommand{\dashlinestretch}{30}
\begin{picture}(6024,489)(0,-10)
\path(462,87)(387,162)
\path(162,387)(237,312)
\path(762,387)(1062,87)
\path(762,87)(837,162)
\path(987,312)(1062,387)
\put(1249.500,237.000){\arc{375.000}{5.3559}{7.2105}}
\put(1774.500,237.000){\arc{375.000}{2.2143}{4.0689}}
\put(2112.000,499.500){\arc{375.000}{0.6435}{2.4981}}
\put(2112.000,-25.500){\arc{375.000}{3.7851}{5.6397}}
\put(2974.500,237.000){\arc{375.000}{2.2143}{4.0689}}
\put(2449.500,237.000){\arc{375.000}{5.3559}{7.2105}}
\path(2637,237)(2787,237)
\put(3312.000,499.000){\arc{374.401}{0.6414}{2.5002}}
\put(3312.000,-25.000){\arc{374.401}{3.7830}{5.6418}}
\path(162,87)(462,387)
\path(3312,162)(3312,312)
\path(5937,462)(6012,462)(6012,12)(5937,12)
\put(3649.500,237.000){\arc{375.000}{5.3559}{7.2105}}
\put(4174.500,237.000){\arc{375.000}{2.2143}{4.0689}}
\put(3837,237){\blacken\ellipse{60}{60}}
\put(3837,237){\ellipse{60}{60}}
\put(4249.500,237.000){\arc{375.000}{5.3559}{7.2105}}
\put(4774.500,237.000){\arc{375.000}{2.2143}{4.0689}}
\put(4587,237){\blacken\ellipse{60}{60}}
\put(4587,237){\ellipse{60}{60}}
\put(5112.000,499.500){\arc{375.000}{0.6435}{2.4981}}
\put(5112.000,-25.500){\arc{375.000}{3.7851}{5.6397}}
\put(5112,312){\blacken\ellipse{60}{60}}
\put(5112,312){\ellipse{60}{60}}
\put(5712.000,499.500){\arc{375.000}{0.6435}{2.4981}}
\put(5712.000,-25.500){\arc{375.000}{3.7851}{5.6397}}
\put(5712,162){\blacken\ellipse{60}{60}}
\put(5712,162){\ellipse{60}{60}}
\path(87,462)(12,462)(12,12)(87,12)
\end{picture}
} } 
  \end{array}
 \]

\end{document}